\documentclass{amsart}

\usepackage[all] {xy}
\usepackage[dvips]{graphicx}
\usepackage{graphics,color}
\usepackage{tikz}
\usepackage[american]{babel}        
\usepackage[T1]{fontenc}          

\usepackage{amsthm,amsmath,mathrsfs}
\usepackage{amssymb,amsthm,amsmath}
\usepackage[dvips]{geometry}
\usepackage{amscd}
\usepackage[latin1]{inputenc}

\usepackage{colortbl}
\usepackage{color,graphics}

\usepackage{psfrag}

\title[ A link between Topological Entropy and Lyapunov Exponents]{A link between Topological Entropy and\\ Lyapunov Exponents}

\author[T. Catalan]{Thiago Catalan}

\theoremstyle{plain}
\newtheorem{mainthm}{Theorem}
\newtheorem{mainclly}[mainthm]{Corollary}
\swapnumbers
\newtheorem{theorem}{Theorem}[section]
\newtheorem{proposition}[theorem]{Proposition}

\newtheorem{lemma}[theorem]{Lemma}

\newtheorem{definition}[theorem]{Definition}

\newcommand{\diff}{\operatorname{Diff}}
\newcommand{\dw}{\operatorname{Diff}^1_{\omega}(M^{2n})}

\newcommand{\Z}{\mathbb{Z}}
\newcommand{\N}{\mathbb{N}}
\newcommand{\R}{\mathbb{R}}

\newcommand{\eps}{\varepsilon}

\begin{document}

\begin{abstract}
We show that a $C^1-$generic non partially hyperbolic symplectic diffeomorphism $f$ has topological entropy equal to the supremum of the sum of the positive Lyapunov exponents of its hyperbolic periodic points. Moreover,  we also prove that $f$ has topological entropy approximated by the topological entropy of $f$ restrict to basic hyperbolic sets. In particular, the topological entropy map is lower semicontinuous in a $C^1-$generic set of symplectic diffeomorphisms far from partial hyperbolicity. 
\end{abstract}

\let\thefootnote\relax\footnote{2000 {\it Mathematics Subject Classification}.
Primary 37J10, 37D30, 37B40, 37C20.}
\let\thefootnote\relax\footnote{{\it Key words and phrases}.
Topological entropy, Lyapunov exponents, homoclinic tangency, generic properties.}
\let\thefootnote\relax\footnote{Work partially supported by CNPq and FAPEMIG}

\maketitle

\section{Introduction}
\label{s.introduction}
An important problem in the theory of Dynamical Systems  is how to measure/compare the complexity of  different dynamics. 
The most famous tool, in this sense, is the topological entropy. This invariant concept (which means it is invariant by conjugacy), roughly,  measures the rate of exponential growth of the number of distinct orbits when we compare  orbits with finite length. Formally, if $f:X\rightarrow X$ is a continuous map defined in a compact metric space $(X,d)$, then we define the {\it topological entropy} of $f$ by:
$$
h_{top}(f)= \lim_{\eps\rightarrow 0} \limsup_{n\rightarrow \infty} \frac{1}{n}\log r(n,\eps),
$$
where  $r(n,\eps)$ is the maximum number of $\eps-$distinct orbits of length $n$. Two points $x$ and $y$ have $\eps-$distinct orbits of length $n$ if there is $0\leq j\leq n$ such that  $d(f^j(x), f^j(y)) > \eps.$

Other way to measure the complexity of dynamics is using Lyapunov exponents. Informally, this number tell us the rate of divergence of nearby trajectories.  Lyapunov exponents measure the asymptotic behavior of dynamics in the tangent space level.  Given a diffeomorphism $f$ over a manifold $M$, we say that
a real number $\lambda(x,f)$ is a {\it Lyapunov exponent} of $x\in M$ for $f$ if there
exists a nonzero vector $v\in T_xM$ such that
$$
\lim_{n\rightarrow \pm\infty} \frac{1}{n}\log \|Df^n(x) \ v\|=\lambda(x,f).
$$

The main result of this paper links topological entropy and Lyapunov exponents in a particular setting.
Namely, we prove that a $C^1-$generic non partially hyperbolic symplectic diffeomorphism $f$ has topological entropy equal to the supremum of the sum of the positive Lyapunov exponents of its hyperbolic periodic points. This was already proved in dimension two by the author together with Tahzibi in  \cite{CT}. Also, it is worth to point out that our result  improves  Theorem B in \cite{CH}, which one gives a lower bound for topological entropy by means of Lyapunov exponents of periodic points. Se also \cite{N}.

Let us make precisely the statement.  

Here we will consider $M^{2n}$ as a compact, connected and boundaryless Riemannian manifold with dimension $2n$. Recall $M^{2n}$ is a symplectic manifold if $M^{2n}$ admits a {\it symplectic form}, i.e., if the manifold has a non-degenerated, closed,  differential $2-$form $\omega$. We denote a symplectic manifold by $(M^{2n},\omega)$. We say a diffeomorphism $f$ acting on $M^{2n}$ is symplectic if it preserves the symplectic form $\omega$, i.e, $f^*\omega=\omega$. The set of $C^1$ symplectic diffeomorphisms is denoted by $\dw$.

Leb $f \colon M^{2n} \to M^{2n}$ be a $C^1-$diffeomorphism and $\Lambda$ be a closed,  $f-$invariant subset of $M$. We say that $\Lambda$ is a {\it partially hyperbolic set}
if there exists a continuous $Df$-invariant splitting
$T_{\Lambda}M=E^s\oplus E^c\oplus E^u$ with non trivial extremal sub-bundles
$E^s$ and $E^u$, such that for every $x\in \Lambda$ and every $m$ large enough:
\begin{itemize}
\item the splitting is \emph{dominated}:
$$
\|Df^m|E^i(x)\| \ \|Df^{-m}| E^j(f^m(x))\|\leq \displaystyle\frac{1}{2},
\ \text{ for any } (i,j)=(s,c),\ (s,u), \ (c,u); \text{ and}
$$
\item the extremal subbundles are hyperbolic:
$$
 \|Df^m|E^s(x)\|\leq \displaystyle\frac{1}{2} \quad \text{ and }\quad
 \|Df^{-m}|E^u(x)\|\leq\displaystyle\frac{1}{2}.
$$
\end{itemize}
If $E^c$ is trivial, then the set $\Lambda$ is called {\it hyperbolic.} When $M$ is a partially hyperbolic set of $f$, we say that {\it $f$ is partially hyperbolic}, and when $M$ is a hyperbolic set of $f$, we say
$f$ is an {\it Anosov diffeomorphism.} We wil denote by $\mathcal{PH}$ the set of $C^1$ partially hyperbolic diffeomorphisms.  It is worth to point out that decompositions satisfying the first item in the definition of partial hyperbolicity is called {\it dominated splitting}.

Now, if $p$ is a periodic point of  a $C^1-$diffemorphism $f$, and if $\lambda$ is a complex eigenvalue of $Df^{\tau(p,f)}$, then $\chi=\frac{1}{\tau(p,f)} \log \|\lambda\|$ is a Lyapunov exponent of $p$ for $f$. In this paper, $\tau(p,f)$ denotes {\it the period} of the periodic point $p$.  We say $\chi$ has multiplicity $m$, if $\lambda$ is an eigenvalue of $Df^{\tau(p,f)}$  with multiplicity $m$. Hence, we define: 
$$
S(p,f)=\sum_{\chi^+_i(p,f)>0} \chi^+_i(p,f),
$$
where the sum is over all positive Lyapunov exponents of $p$ for $f$, counting multiplicity. Using this notation, we also define
$$
S(f)=\sup\{S(p,f); \text{ where } p \text{ is a hyperbolic periodic point of } f\}. 
$$
Recall we say a {\it periodic point is hyperbolic} if  the orbit of $p$ is a hyperbolic set of $f$.  Since these points are robust (they are preserved after perturbations of the diffeomorphism), it is possible to conclude that the functional $S(.)$ is a lower semicontinuous map for $C^1-$diffeomorphisms.

Finally, we can state the main result of the paper:

\begin{mainthm}
\label{entropy}
There exists a residual subset $\mathcal{R}\subset
\dw$ such that if $f\in\mathcal{R}$ is not partially hyperbolic, then
$$
h_{top}(f)= S(f).
$$
\end{mainthm}

It is important to remark that this equivalence  does not happen for Anosov diffeomorphisms, in general. In fact, there is an open and dense subset of Anosov symplectic diffeomorphisms $f$ where $h_{top}(f)<S(f)$. Hence, we put a question: What happens for partially hyperbolic diffeomorphisms with non trivial center bundle?  

Another important problem in smooth theory is about the regularity of the topological entropy. Using Theorem \ref{entropy} will be possible to prove something in this direction.  For instance, Newhouse has proved in \cite{Ne} that topological entropy map varies upper semi continuously for $C^{\infty}$ diffeomorphisms ($C^{\infty}$ topology), which is not true in the $C^1-$topology. See Theorem C in \cite{CT}. On the other hand, Katok shows that topological entropy varies lower semi continuously for $C^2-$diffeomorphisms ($C^2-$topology) over surfaces. See \cite{KH} for a proof. In particular, the topological entropy is continuous for $C^{\infty}$-surface diffeomorphisms ($C^{\infty}$ topology). To prove the lower semi continuity, Katok proved that the topological entropy of a $C^2-$diffeomorphism $f$ is approximated by topological entropy of $f$ restrict to hyperbolic sets. In the lack of partial  hyperbolicity we can find a similar result for $C^1-$generic symplectic diffeomorphisms. 

Given a diffeomorphism $f$ over $M^{2n}$, we define:

$$
H(f)=\sup\{h_{top}(f|\Lambda); \, \text{ where } \Lambda \text{ is a basic hyperbolic set of } f\}.
$$
By {\it basic hyperbolic set } we mean a hyperbolic set which is isolated, i.e., there is a neighborhood $U$ of $\Lambda$ such that $\Lambda=\cap_{j\in \Z} f^j(U)$.

\begin{mainthm}
There is a residual subset $\mathcal{R}\subset \dw$, such that if $f\in\mathcal{R}$ is not partially hyperbolic, then 
$$
h_{top}(f)=H(f).
$$
\label{t.katok}\end{mainthm}

Therefore, as Katok, by invariance of the topological entropy and robustness of basic hyperbolic sets,  the previous theorem implies the existence of points where the topological entropy map has a regularity. More precisely: 

\begin{mainclly}
There is a non empty residual subset $\mathcal{R}\subset \dw-cl(\mathcal{PH})$, such that if $f\in\mathcal{R}$, then the topological entropy map, $h_{top}: \diff^1(M)\rightarrow \R^+$, is lower semi-continuous at $f$.
\label{cor.}\end{mainclly}

Recall that the topological entropy map is locally constant for Anosov diffeomorphism. However, we do not know a similar result to Corollary \ref{cor.} for partially hyperbolic diffeomorphism with non trivial center bundle.

This paper is organized in the following way:   in section 2, we prove Theorem \ref{entropy}  using Propositions \ref{diag.equal}, \ref{fulltangency} and \ref{tech.prop}. In section 3, we prove Theorem \ref{t.katok} and Corollary \ref{cor.}. Finally,  in section 4, 5 and 6 we prove Propositions  \ref{diag.equal}, \ref{fulltangency} and \ref{tech.prop}, respectively.

\vspace{0,3cm}
{\bf Acknowledgements:} I would like to thanks A. Arbieto, S. Crovisier, V. Horita, and R. Potrie for useful conversations and constant encouragement. To IM-UFRJ and the Abdus Salam ICTP for the kind hospitality during preparation of this work. This work was partially supported by CNPq and FAPEMIG.

\section{Proof of Theorem \ref{entropy}}

We divide the proof of Theorem \ref{entropy} in two results. The first one shows that the functional $S(.)$ is generically an upper bound for the topological entropy. And in the second one we prove that generically any symplectic non partially hyperbolic diffeomorphism $f$ has topological entropy lower bounded by $S(f)$.  More precisely, Theorem \ref{entropy} follows directly from the following two results:

\begin{theorem}
There exists a residual subset $\mathcal{R}\subset
\dw$ such that if $f\in\mathcal{R}$, then
$$
h_{top}(f)\leq S(f).
$$
\label{hleqs}\end{theorem}

\begin{theorem}
There exists a residual subset $\mathcal{R}\subset
\dw$ such that if $f\in\mathcal{R}$ is not partially hyperbolic, then
$$
h_{top}(f)\geq S(f).
$$
\label{hgeqs}\end{theorem}

Before we start the proofs we recall some perturbative results. First, a symplectic version of Franks' lemma \cite{F}.

\begin{lemma}[Lemma 5.1 (Symplectic Franks' lemma) in \cite{HT})]. Let $f\in \dw$ and $E$ a finite $f$-invariant set. Assume that B is a small symplectic perturbation of $Df$ along $E$. Then for every neighborhood $V\subset M$ of $E$ there is a symplectic diffeomorphism $g$ arbitrarily $C^1$-close to $f$ coinciding with $f$ on $E$ and out of $V$, and such that $Dg$ is equal to $B$
on $E.$
\label{franks}\end{lemma}

And a pasting lemma of Arbieto and Matheus \cite{AM}. 

\begin{theorem}[Pasting lemma] 
If $f$  is a $C^k-$symplectic  diffeomorphism over $M$ $(k\geq 1)$, and $x\in M$, then for every $\eps>0$ there exists a $C^1$ symplectic diffeomorphism $g$ $\eps-C^1$ close to $f$ such that for small neighborhoods of $x$, $V\subset U$, $g|U^c=f$ and $g|V=Df(x)$ (in local coordinates). 
\label{pasting lemma}\end{theorem}

\subsection{Proof of Theorem \ref{hleqs}} 
We remark that Theorem \ref{hleqs} is a consequence of Ruelle's inequality and a symplectic version of a result of Abdenur, Bonatti and Crovisier (Theorem 3.8 \cite{ABC}). 

Given an ergodic invariant measure $\mu$ of a diffeomorphism $f$ on $M^{2n}$, $\mu\in \mathcal{M}_e(f)$,  there are real numbers (not necessary distinct)
$\chi_1\leq \chi_2\leq \ldots \leq \chi_{2n}, $ which are the Lyapunov exponents for $\mu-$almost every point $x\in M$. We define the {\it Lyapunov vector} of $\mu$ by: $L(\mu)=(\chi_1,\ldots,\chi_{2n})$. Also, given a periodic point $p$ of a diffemorphism $f$, we denote by {\it $\mu_p$ the Dirac measure defined in the orbit of $p$.}

\begin{theorem}[Symplectic version of Theorem 3.8 \cite{ABC}]
There exists a residual subset $\mathcal{R}\subset\dw$ such that if $f\in \mathcal{R}$ then 
 for an ergodic measure $\mu$ of $f$, there is a sequence of periodic points $p_n$ of $f$ such that $\mu_{p_n}$ converge to $\mu$ in the weak topology, and moreover the Lyapunov vectors $L(\mu_{p_n})$ also converge to the Lyapunov vector $L(\mu)$.  
\label{t.abc}\end{theorem}

{\it Proof:}
The original proof given in \cite{ABC} also works here, since there exist symplectic versions of the perturbative results used there.  

$\hfill\square$

To prove Theorem \ref{hleqs}, we also need the following technical lemma: 

\begin{lemma}
There exists a residual subset $\mathcal{R}\subset\dw$ such that if $f\in \mathcal{R}$ then 
$$
S(f)= \sup_{p \in Per(f)} S(p,f). 
$$
\label{lemma.S.s}\end{lemma}

{\it Proof:} By Kupka-Smale Theorem, generically a symplectic diffeomorphism has only hyperbolic or m-elliptic periodic points. Recall that $p$ is a $m-$ellipitic periodic point, $0<\leq m \leq n$, if $Df^{\tau(p,f)}$ has exactly $2m$ non-real and simple eigenvalues of modulus one, and all other eigenvalues has modulus different from one. These periodic points are also robust. By Kupka-Smale Theorem  we have, in fact, that if we fix a period, there exists  an open and dense subset of symplectic diffeomorphisms for which there are only finite periodic points (hyperbolic or m-elliptic) of such period.  Hence, by standard generic arguments, the proof of the lemma follows naturally from the following claim: 

\vspace{0,2cm}
{\it Claim: For any $\eps>0$ and any $m-$elliptic periodic point $p$ of a symplectic diffeomorphism $f$, there exists a symplectic diffeomorphisms $\tilde{f}$ arbitrary close to $f$, and a hyperbolic periodic point $\tilde{p}$ of $\tilde{f}$ such that $L(\tilde{p})$ is $\eps-$close to $L(p)$.}

\vspace{0,1cm}
In fact, after a perturbation of $f$, we can suppose that the argument of all non-real eigenvalues with modulus equal to one of $Df^{\tau(p)}(p,f)$ are rational. Now, using pasting lemma, Lemma \ref{pasting lemma}, to perform  a local perturbation of $f$, we can also suppose $f^{\tau(p,f)}=Df^{\tau(p,f)}(p)$ in a small neighborhood of $p$ (in local coordinates). Hence, we can conclude the existence of a non-hyperbolic periodic point $\tilde{p}$ close to $p$ having only real eigenvalues, such that $Df^{\tau(\tilde{p},f)}(\tilde{p})=(Df^{\tau(p),f}(p))^k$,  for some $k\in \N$. In particular the Lyapunov vector of $\tilde{p}$ is the same of the Lyapunov vector of $p$. Moreover, we can choose  linear maps $A_i$ arbitrary close to $Df^i(\tilde{p})$, $1\leq i\leq \tau(\tilde{p},f)$, such that $\prod A_i$ has no eigenvalue with absolute value equal to one.  Thus, we can use Franks Lemma, to find $\tilde{f}$ arbitrary close to $f$ such that  $\tilde{p}$ still is a periodic point of $\tilde{f}$, and $D\tilde{f}^{i}(\tilde{p})=A_i$, for $1\leq i\leq \tau(\tilde{p},\tilde{f})$, which implies $\tilde{p}$ is hyperbolic for $\tilde{f}$. By continuity of the Lyapunov vector the proof of the claim is finished, and thus the proof of the lemma.  

$\hfill\square$

{\it Proof of Theorem \ref{hleqs}:}

Let $\mathcal{R}$ be a residual subset for which is true Theorem \ref{t.abc} and Lemma \ref{lemma.S.s}. Consider now $f\in \mathcal{R}$. Then, for any $\mu\in \mathcal{M}(f)$ there exists a sequence of periodic points $p^{\mu}_{ n}$ such that the dirac measures on the orbit of $p^{\mu}_n$ converge to $\mu$ and the Lyapunov
exponents of $p^{\mu}_n$ converge to the Lyapunov exponents of $\mu.$ By variational principle $h_{top} (f) = \sup_{\mu} h_{\mu}(f)$, where the supremum is
over all f-invariant ergodic probability measures. By Ruelle's inequality:   $h_{\mu}(f) \leq \sum \chi_i^+(\mu)$, where the sum is over all positive
Lyapunov exponents of $\mu$, counting multiplicity.  Then we have that
$$
h_{top}(f) \leq \sup_{p \in Per(f)} S(p,f).
$$
Since $f\in \mathcal{R}$, by Lemma \ref{lemma.S.s}, we have $h_{top}(f)\leq S(f)$.

$\hfill\square$

\subsection{Proof of Theorem \ref{hgeqs}}

 Theorem \ref{hgeqs}  will be a consequence of some propositions. We prefer to conclude the proof of the theorem and postpone the proofs of the propositions for other sections.

\vspace{0,2cm}

{\it Proof of Theorem \ref{hgeqs}:}

Given a positive integer $k>0$,  we define
$\mathcal{B}_{k}\subset \dw-cl(\mathcal{PH})$ as the subset of diffeomorphisms
$g$ such that
$$
h_{top}(g)>S(g)-\frac{1}{k}.
$$
The proof of the theorem follows directly of the following claim, by considering: 
$$
\mathcal{R}=\bigcap_{k\in \N} (\mathcal{B}_k\cup \mathcal{PH}).
$$

\vspace{0,2cm}
{\it Claim: For every $k\in \N$,  $\mathcal{B}_{k}$ is an open and dense subset in $\dw-cl(\mathcal{PH})$.  }

\vspace{0,1cm}
{\it Proof of the Claim:}

First, since $S(\cdot)$ is a lower semicontinuous map defined in $\dw$,
there is a residual subset  $\mathcal{R}\subset \dw$ where
$S(\cdot)$ is continuous. 

Now, let $f\in \mathcal{R}\cap (\dw-cl(\mathcal{PH}))$, and $\mathcal{V}\subset   \dw-cl(\mathcal{PH})$ be a small neighborhood of $f$, such that
\begin{equation}
S(f)>S(g)-\frac{1}{5k}, \text{ for every } g\in \mathcal{V}.
\label{t.eq0}\end{equation}

By definition of $S(f)$ there exists a hyperbolic periodic point $p$ of $f$ such that 
\begin{equation}
S(f)>S(p,f)-\frac{1}{5k}. 
\label{t.eq1}\end{equation}

We need now the following result: 

\begin{proposition} 
Let $f\in \dw$ be not approximated by  partially hyperbolic symplectic diffeomorphisms and $p$ be a hyperbolic periodic point of $f$. For any $C^1-$neighborhood $\mathcal{U} \subset \dw$ of $f$ and $\eps>0$, there exists a diffeomorphism $g\in \mathcal{U}$ and a hyperbolic periodic point $\tilde{p}$ of $g$ such that: 
$$
|S(p,f)-n \chi_{min}^+(\tilde{p},g)|<\eps.
$$
Where $\chi_{min}^+(\tilde{p},g)$ is the smallest positive Lyapunov exponent of $\tilde{p}$ for $g$.
\label{diag.equal}\end{proposition}

Hence, by Proposition \ref{diag.equal}, we can perturb $f$ to $f_1\in \mathcal{V}$ such that there is a hyperbolic periodic point $\tilde{p}$ of $f_1$ satisfying: 
\begin{equation}
|S(p,f)-n \chi_{min}^+(\tilde{p},f_1)|<\frac{1}{5k}.
\label{t.eq2}\end{equation}

Also, since $f_1$ still is not approximated by partially hyperbolic diffeomorphisms, it is possible to find a full homoclinic tangency for $\tilde{p}$ after perturbing $f_1$. See the next result. We say that a periodic point  $\tilde{p}$ of a symplectic diffeomorphism $f_1$ {\it exhibits a full homoclinic tangency} if there exists a point $q$ ({\it homoclinic point}) which is a point of intersection between stable and unstable manifold of $\tilde{p}$, $W^s(\tilde{p},f_1)$ and $W^u(\tilde{p},f_1)$, respectively, such that $T_qW^s(\tilde{p})=T_qW^u(\tilde{p})$.  Recall, that hyperbolic periodic points of symplectic diffeomorphisms over $2n-$manifolds have index equal to $n$, i.e., the dimension of the stable manifold is $n$.

\begin{proposition}
Let $f\in \diff^1_{\omega}(M^{2n})$ be a symplectic diffeomorphism which is not $C^1$ approximated by  partially hyperbolic symplectic diffeomorphisms. If $p$ is a hyperbolic periodic point of $f$ then for any $C^1$ small neighborhood $\mathcal{U}\subset \dw$ of $f$, there exists a diffeomorphism $g\in \mathcal{U}$ exhibiting a full homoclinic tangency for $p(g)$. Where $p(g)$ is the analytic continuation of $p$ for $g$. 
\label{fulltangency}\end{proposition}

Thus, using Proposition \ref{fulltangency}, we can find $f_2\in \mathcal{V}$ arbitrary close to $f_1$ such that $\tilde{p}(f_2)$ has a full homoclinic tangency for $f_2$. For convenience, we still denote $\tilde{p}(f_2)$ by $\tilde{p}$. 

The following result is the key technical point is this proof.

\begin{proposition}
Let $f\in \diff^1_{\omega}(M^{2n})$  and let $p$ be a hyperbolic periodic point  of $f$. 
If $f$ exhibits  a full homoclinic tangency $q$ for $p$,   then for any neighborhood $\mathcal{U}$ of $f$ and any $\eps>0$, there exists a diffeomorphism $g\in \mathcal{U}$ such that $p$ still is a hyperbolic periodic point of $g$ and moreover there exists a basic hyperbolic set $\Lambda_p$ of $g$ such that 
$$
h(g|\Lambda_p)\geq n\ \chi^+_{min}(p,g)-\eps. 
$$     
\label{tech.prop}\end{proposition}

By Proposition \ref{tech.prop}, there exists a diffemorphism $f_3\in  \mathcal{V}$ arbitrary close to $f_2$, and a hyperbolic set $\Lambda_{\tilde{p}}$ of $f_3$ such that

\begin{equation}
h(f_3|\Lambda_{\tilde{p}})\geq n\ \chi^+_{min}(\tilde{p},f_3)-\frac{1}{5k}. 
\label{t.eq3}\end{equation}

Therefore, if $\mathcal{U}\subset \mathcal{V}$ is a small enough neighborhood of $f_3$, then for any $g\in \mathcal{U}$ 
there is a continuation of the hyperbolic basic set $\Lambda_{\tilde{p}}$ which
we denote by $\Lambda(\tilde{p},g)$. Thus, using respectively,  properties of entropy, equation (\ref{t.eq3}),  continuity of the Lyapounov exponents, equation (\ref{t.eq2}), equation (\ref{t.eq1},), and equation (\ref{t.eq0}), we have
\begin{align}\label{t.eq.4}
 h_{top}(g) \geq h_{top}(g|\Lambda(\tilde{p}, g))
&= h_{top}(f_3|\Lambda_{\tilde{p}})
\\
\nonumber & \geq n\ \chi^+_{min}(\tilde{p},f_3)-\frac{1}{5k}.
\\
\nonumber &\geq n\ \chi^+_{min}(\tilde{p},f_1)-\frac{2}{5k}
\\
\nonumber &\geq S(p,f)-\frac{3}{5k}
\\
\nonumber &\geq S(f)-\frac{4}{5k}
\\
\nonumber &\geq S(g)-\frac{1}{k}. 
\end{align}

Hence, every $g\in \mathcal{U}$  belongs to $\mathcal{B}_{k}$. Since $\mathcal{R}$ is dense in $\dw-cl(\mathcal{PH})$, we can conclude that $\mathcal{B}_{k}$ contains an open and dense subset in
 $\dw-cl(\mathcal{PH})$. 
 
 $\hfill\square$

\section{Regularity of topological entropy}

First we prove Theorem \ref{t.katok}.

\vspace{0,2cm}
{\it Proof of Theorem \ref{t.katok}:}
Clearly, $h_{top}(f)\geq H(f)$. Let we prove now that $H(f)\geq h_{top}(f)$  for symplectic generic diffeomorphisms  in the lack or partial hyperbolicity. For that, we define  $\mathcal{C}_k$ the set formed by symplectic difeomorphisms such that $H(g)\geq S(g)-1/k$, for any $k\in \N$. Looking to equation (\ref{t.eq.4}), we can conclude analogously we did in the proof of Theorem \ref{hgeqs},  that $\mathcal{C}_k\cap (\dw-cl(\mathcal{PH}))$ contains an open and dense subset. Hence, 

$$
\mathcal{R}_1=\bigcap_{k\in \N} \mathcal{C}_k\cup \mathcal{PH} \text{ is  a residual subset of } \dw.
$$ 
Let $\mathcal{R}_2\subset \dw$ be the residual subset  given by Theorem \ref{entropy}.  Now, if $f\in \mathcal{R}=\mathcal{R}_1\cap \mathcal{R}_2$ is not partially hyperbolic, then $h_{top}(f)= S(f)$ and $H(f)\geq S(f)$, by choice of $\mathcal{R}_1$ and $\mathcal{R}_2$, respectively. Which implies that $H(f)\geq h_{top}(f)$. And thus, $H(f)=h_{top}(f)$ for any $f\in\mathcal{R}$.

$\hfill\square$

\vspace{0,2cm}

{\it Proof of Corollary \ref{cor.}:}
Let $\mathcal{R}$ be the residual subset given by Theorem \ref{t.katok} and let $f\in\mathcal{R}$.  By definition of $H(f)$, given any $\eps>0$ there exists a basic hyperbolic set $\Lambda$ of $f$ such that $h_{top}(f|\Lambda)>h_{top}(f)-\eps$. By continuity of the hyperbolicity, for any $C^1-$diffeomorphism $g$, not necessary symplectic, close to $f$ there is a basic hyperbolic set $\Lambda(g)$ which is close and topologically conjugated to $\Lambda$. Thus we have: 

\begin{align*}
h_{top}(g)&\geq h_{top}(g|\Lambda(g))
\\
&=h_{top}(f|\Lambda)
\\
&\geq h_{top}(f)-\eps,
\end{align*}
for any $g$, $C^1-$close to $f$. That is, $f$ is point of lower semicontinuity of the topological entropy map in $\diff^1(M)$.

$\hfill\square$

\section{Proof of Proposition \ref{diag.equal}}

This proposition is a consequence of some results about {\it Periodic linear systems with transitions}. This kind of system was introduced by Bonatti, Diaz and Pujals in \cite{BDP}. Many of the results in \cite{BDP} is also true in the symplectic setting, proved by Horita and Tahzibi in \cite{HT}. In this section, we will show how to use their results and arguments to prove Proposition \ref{diag.equal}. We also use some arguments developed in  \cite{CH}. 

We remark that since  the theory of Periodic linear systems is well known nowadays, we will give only a sketch of the proof. 

Let we recall some definitions and some basic results used in this section. 

Given a diffemorphism $f$ over $M^{2n}$, we say that two hyperbolic periodic points $p_1$ and $p_2$ are {\it homoclinically related} if  $W^s(p_1,f)$ and $W^u(p_2,f)$ intersects transversally $W^u(p_1,f)$ and $W^s(p_2,f)$, respectively. That is, the stable and unstable manifolds of $p_1$ and $p_2$ have non-empty transversally intersections. We denote by $H(p_1,f)$ the closure of all hyperbolic periodic points of $f$ homoclinically related to $p_1$. This set is called by {\it homoclinic class of} $p_1$. It is worth to point out that $f$ restrict to a homoclinic class is {\it topological transitive}, i.e., this set contains a forward dense orbit of $f$. 

The following result due to M-C. Arnaud, Bonatti and Crovisier says that any homoclinic class of a generic symplectic diffemorphism coincides with the whole manifold ambient. This is a symplectic version of the main result in \cite{BC}.

\begin{theorem}[Theorem 1 in \cite{ArBC}]
There is a residual subset $\mathcal{R}\subset \dw$ such that if $f\in \mathcal{R}$ and $p$ is a hyperbolic periodic point of $f$, then $H(p,f)=M$. 
\label{ArBC}\end{theorem}

\subsection{Periodic Linear systems with symplectic transitions}
In this subsection we  will recall some basic definitions. A  {\it periodic symplectic linear system (or linear symplectic cocycle over $f$)} is a 4-tuple $\mathcal{P}=(\Sigma, f, \mathcal{E},A)$, where $f$ is a homeomorphism,
$\Sigma$ is an infinite set of periodic points of $f$, $\mathcal{E}$ is a symplectic vector bundle defined over $\Sigma$ (i.e., each fiber $\mathcal{E}(x)$ is a symplectic vector space), and $A:\mathcal{E}\rightarrow \mathcal{E}$ is   such that $A(x):\mathcal{E}_x\rightarrow \mathcal{E}_{f(x)}$ is a symplectic linear isomorphism for each $x$.
Also, a system $(\Sigma, f, \mathcal{E},A)$ can  be considered as a system of matrices $(\Sigma, f, A)$. In particular, we denote by
$$
M_A (x) = A(f^{n-1}(x)) A(f^{n-2}(x)) \ \dots \   A(x), $$
if $x\in \Sigma$ and $n$ is the period of  $x$.
A periodic linear system is {\it diagonalizable} at the point $x\in \Sigma$ if  $M_A(x)$
has only real eigenvalues of multiplicity one.

Given a set $\mathcal{B}$, a {\em word} with letters in $\mathcal{B}$
is a finite sequence of elements of $\mathcal{B}$.
The product of the word $[a]= (a_1, \dots , a_n)$ by
$[b]= (b_1, \dots , b_m)$ is the word $(a_1, \dots , a_n, b_1, \dots, b_m)$.
We say a word is {\em not a power} if $[a] \neq [b]^k$ for every word
$[b]$ and $k > 1,$ and we say that two words $[a]$ and $[b]$ are $\eps-$close if they have the same length, and the correspondents letters are $\eps-$close.  With this notation, for a periodic symplectic linear system $(\Sigma, f, \mathcal{E}, A)$ if we consider the word
$[M_A (x)] = (A(f^{n-1}(x)), \dots , A(x)),$
where $n$ is the period of $x\in \Sigma$, then
the matrix $M_A(x)$ is the product of the letters of the word $[M_A (x)]$.

\begin{definition}[Definition 1.6 of \cite{BDP}]
\label{d.transition} Given $\varepsilon > 0$, a periodic symplectic linear system $(\Sigma, f, \mathcal{E},A)$ {\em admits $\varepsilon$-transitions} if for every finite family of points $x_1, \dots , x_n = x_1 \in \Sigma$, 
and for any $(i, j) \in \{1,\dots
,n\}^2$ there exist $k(i, j) \in \mathbb{N}$ and a finite word
$[t^{i,j}]= (t_1^{i,j}, \dots , t_{k(i,j)}^{i,j})$ of symplectic matrices, satisfying the following properties:
\begin{enumerate}
\item For every $m \in \mathbb{N}$, $\imath = (i_1, \dots, i_m) \in \{1,\dots,n\}^m$,
and $\alpha = (\alpha_1,\dots ,\alpha_m) \in \mathbb{N}^m$
consider the word
\begin{align*}
[W(\imath, \alpha)] & = [t^{i_1,i_m}][M_A(x_{i_m})]^{\alpha_m}
[t^{i_m,i_{m-1}}][M_A(x_{i_{m-1}})]^{\alpha_{m-1}} \dots \\
& \quad \dots [t^{i_2,i_1}][M_A(x_{i_1})]^{\alpha_1},
\end{align*}
where the word
$[W(\imath, \alpha)]$ is not a power.
Then there is $x(\imath, \alpha) \in \Sigma$ such that
\begin{itemize}
\item The length of $[W(\imath, \alpha)]$ is the period of
$x(\imath,\alpha)$.
\item The word $[M_A (x(\imath, \alpha))]$ is $\varepsilon$-close to
$[W(\imath, \alpha)]$ and there is an $\varepsilon$-symplectic perturbation
$\tilde A$ of $A$ such that the word $[M_{\tilde{A}} (x(\imath,\alpha))]$
is $[W(\imath, \alpha)]$.
\end{itemize}

\item One can choose $x(\imath, \alpha)$ such that the distance between
the orbit of $x(\imath, \alpha)$ and any point $x_{i_k}$ is bounded by some
function of $\alpha_k$ which tends to zero as $\alpha_k$ goes to infinity.
\end{enumerate}
\end{definition}

In our setting, results about periodic linear systems are useful because of the following result, 
which is a symplectic version of [Lemma 1.9 in [BDP]].

\begin{lemma}[Lemma 4.5 in \cite{HT}]  Let $f$ be a symplectic diffeomorphism and let $p$ be a hyperbolic periodic point of $f$. The derivative $Df$ induces a continuous periodic symplectic linear system with transitions on the set $\Sigma$ formed by the hyperbolic periodic points homoclinically related  to $p$. 
\label{ex}\end{lemma}

A nice property of periodic linear systems with transitions is the fact that we have a diagonalizable arbitrarily small perturbation of $A$ defined on a dense subset of $\Sigma$, as we can see in Lemma 4.16 in \cite{BDP} (see Lemma 4.7 in \cite{HT} for a symplectic version). Recall we say that a  point $x\in \Sigma$ is {\it diagonalizable for some linear periodic system}, if the matrix $M_A(x)$ is diagonalizable. However, their proof could be improved. That is, we can proceed in the same way as in the proof of Lemma 4.3 in \cite{CH} (see also Remark 4.4 in \cite{CH}), to prove that these diagonalizable perturbations keeps, in some sense, the Lyapunov exponents. More precisely, it is true the following result:

\begin{lemma} Let $(\Sigma, f, \mathcal{E},A)$ be a periodic symplectic linear system. For any $\eps>0$, and $x\in \Sigma$ there exists  $y\in \Sigma$ and  a symplectic perturbation $\tilde{A}$ of $A$ defined on the orbit of $y$, such that $M_{\tilde{A}}(y)$ is diagonalizable and moreover: 
$$
\left|S(x,A)-S(y,\tilde{A})\right|<\eps.
$$
\label{exp}\end{lemma}
We refer the reader to \cite{CH} to see details about the proof of the previous lemma.

Now, given a diagonalizable periodic point $x$ of a periodic symplectic linear system $(\Sigma, f, \mathcal{E}, A)$, let $\lambda_1\leq \ldots\leq\lambda_{2n}$ be the $2n-$distinct eigenvalues of $M_A(x)$. Hence, we will consider the following decomposition $\mathcal{E}(x)=E_1\oplus \ldots \oplus E_{2n}(x)$, where $E_i$ is the one-dimensional eigenspace with respect to $\lambda_i$, for $1\leq i\leq 2n$. 

We will see now, how to use the existence of specific periodic points, namely periodic points having complex eigenvalue of rank $(i, \ i+1)$, to mixing the Lyapunov exponents of a diagonalizable periodic point $x$.

Before, we point out that as done for diffeomorphisms, we can define similar dominated splittings  and partial hyperbolicity on $\mathcal{E}$ for periodic linear systems $(\Sigma, f, \mathcal{E}, A)$. 

\begin{definition}[see Definition 2.2 of \cite{BDP}] We say a periodic linear system $(\Sigma, f, \mathcal{E}, A)$
has {\it a complex eigenvalue of rank $(i,\ i+1)$} if there is $x\in \Sigma$ such that the matriz $M_A(x)$ has a complex eigenvalue of rank $(i,\ i+1)$. That is, if there is a $M_A(x)-$invariant dominated splitting of $\R^{2n}$, $F\oplus G\oplus H$, such that: 
\begin{itemize}
\item every eigenvalue $\sigma$ of $M_A(x)|F$ (resp. $M_A(x)|H$) has modulus $|\sigma|<|\lambda|$ (resp. $|\sigma|>|\lambda|$),

\item $dim(F)=i-1$ and $dim(H)=2n-i-1$, 

\item the plane $G$ is the eigenspace of $\lambda$. 
\end{itemize}
\end{definition}

The following result is a new statement of Proposition 7.1 in \cite{HT}: 

\begin{proposition}
Let $(\Sigma, f, \mathcal{E}, A)$ be a periodic symplectic linear system with symplectic transitions. Fix $\eps_0>0$ and assume that a symplectic $\eps_0-$perturbation of $A$ has a complex eigenvalue of rank $(i, i+1)$ for some $i\in\{1,\ldots, 2N-1\}$. Then for every $0<\eps_1<\eps_0$, if $x\in \Sigma$ is a diagonalizable periodic point, then there is a $\eps_1-$transition $[t^i]$ from $x$ to itself such that: 

\begin{itemize}
\item There is a symplectic ($\eps_0+\eps_1-$)perturbation $[\tilde{t}^i]$ of the transition $[t^i]$ such that the corresponding matrix $\tilde{T}^i$ satisfies
\begin{itemize}
\item[-] $\tilde{T}^i(E_j(x))=E_j(x)$ if $j\not\in [i,\ i+1, \ i^*, (i+1)^*]$, 

\item[-] $\tilde{T}^i(E_i(x))=E_{i+1}(x)$ and $\tilde{T}^{i}(E_{i+1}(x))=E_{i}(x)$ 
\end{itemize}
\end{itemize}
\label{HT.bla}\end{proposition}

Using this proposition we can proceed in the same way as in the proof of Proposition 5.6 in \cite{HT} to conclude:

\begin{proposition}
Let $(\Sigma, f, \mathcal{E}, A)$ be a periodic symplectic linear systems with symplectic transitions. Given $
\eps>0$, if there exists $\delta>0$ such that for any $i\in\{1,\ldots, 2n-1\}$ there is a symplectic $\delta-$perturbation of $A$ having a complex eigenvalue of rank $(i, i+1)$,  then for every diagonalizable point $x\in \Sigma$ there is a point $y\in \Sigma$ and a $\eps-$perturbation $\tilde{A}$ of $A$ defined in the orbit of $y$, such that
$$\left|n\chi^+_{min}(y,\tilde{A})-S(x,A)\right|<\eps.$$ 
\label{blableblu}\end{proposition}

Since the previous proposition is the main technical result in this subsection, we decide to put a proof. 

\vspace{0,2cm}
{\it Proof of Proposition \ref{blableblu}: } Fix $\delta=\eps/3$. Considering $\eps_0=\delta$ and $\eps_1=\eps/3$, let $[t^i]$ be the $\eps_1-$transition from $x$ to itself and  $[\tilde{t}^i]$ the $(\eps_0+\eps_1)-$perturbation of $[t^i]$  given by Proposition \ref{HT.bla}, for any $0<i<2n$. For convenience, we suppose $[t^0]=[t^{2n}]=Id$. 

Note that the maps $[\tilde{t}^i]$ are transpositions acting on the finite set $\{E_j\}_{1\leq j \leq 2n}$ of eigenspaces of $M_A(x)$ which interchanges $E_i$ and $E_{i+1}$. Hence, as any permutation is a composition of transpositions, for every $0\leq k\leq n$ there exists a transition map $[S_k]$, obtained by composing some $[t^i]'$s, such that if we consider $[\tilde{S}_k]$ as the similar word formed by the perturbations $[\tilde{t}^i]$ instead of $[t^i]$, then $\tilde{S}_k(E_{j+n}) = E_{(j+k)+n}$, for any $1\leq j\leq n$, and considering $(j+k)$  in the cyclic group $\Z/n\Z$. Also, recall $\tilde{S}_k$ is the matrix corresponding to $[\tilde{S}_k]$. 

Now, we consider $[W_{k,m}]=[S_{n-k}][M_A(x)]^m[S_k]$, and $[\tilde{W}_{k,m}]=[\tilde{S}_{n-k}][M_A(x)]^m[\tilde{S}_k]$, for any $0\leq k\leq n$. And using such words, we define:

 $$[M_m]= [W_{n-1,m}]\circ\ldots \circ [W_{1,m}]\circ [W_{0,m}]\; \text{ and } \; [\tilde{M}_m]= [\tilde{W}_{n-1,m}]\circ\ldots \circ [\tilde{W}_{1,m}]\circ [\tilde{W}_{0,m}]. $$ 

By definition of periodic linear systems with transitions, there exists a point $y_m\in \Sigma$ such that $[M_A(y_m)]$ is $\eps/3-$close to $[M_m]$, in particular   $\eps-$close to $\tilde{M}_m$. Note that $\tilde{S}_{n-k}\tilde{S}_k$ and $\tilde{M}_m$  keeps invariant any $E_j$, $n< j\leq 2n$. We denote by $\mu_{j,k}$ and  $\tilde{\lambda}_j$ the eingenvalues of $\tilde{S}_{n-k}\tilde{S}_k|E_j$ and $\tilde{M}_m|E_j$, respectively. With this notation we have the following: 
 
 $$\tilde{\lambda}_j=\sum_{j=n+1}^{2n} \lambda_{j}^{m}+ \sum_{k=0}^{n-1} \mu_{j,k}, \ \text{ for any } n+1\leq j\leq 2n.$$

Note now, that the period of $y$ is: $\tau(y)=nm\tau(x)+K$, where $K$ is a constant which depends only on the length of the words $[S_k]$. Thus, if $\tilde{A}_m$ is defined in the orbit of $y_m$ such that $[M_{\tilde{A}}]=[\tilde{M}_m]$, we have that 

$$\chi_j(y_m,\tilde{A}_m)=\displaystyle\frac{1}{n}\sum_{i=n+1}^{2n} \chi_i(x,A)+\frac{1}{\tau(y_m)}\sum_{k=0}^{n-1}\log \mu_{j,k},\ \text{ for any } n+1\leq j\leq 2n.$$

Since $\tau(y_m)$ goes to infinity when $m$ goes, the proof is finished taking $m$ large enough, and considering $y=y_m$ and $\tilde{A}=\tilde{A}_m$. 

$\hfill\square$

Finally, to finish this subsection we state a result which is a directly consequence of Proposition 5.5 in \cite{HT}. 

\begin{proposition}  Let $(\Sigma, \ f, \ \mathcal{E}, \ A)$  be a continuous periodic 2n-dimensional linear system with symplectic transitions such that its norm $\|A\|$ is bounded by positive integer $K$. Given $\eps>0$, If $(\Sigma, \ f, \ \mathcal{E}, \ A)$ admits no dominated splitting, then for any  $1\leq i\leq 2n-1$  there is a $\eps-$perturbation $\tilde{A}$ of $A$ having a complex eigenvalue of rank $(i,\ i + 1)$. 
\label{every.rank}\end{proposition}

\subsection{Proof of Proposition \ref{diag.equal}:}

Let $\eps>0$ be arbitrary small  and $\mathcal{U}$ be a small enough neighborhood of $f$ such that any diffeomorphism inside $\mathcal{U}$ is no partially hyperbolic. In particular, every diffeomorphism in this neighborhood admits no dominated splitting. Recall that for symplectic diffeomorphisms,  partial hyperbolicity  is equivalent to the existence of a dominated splitting on $M$.  After a perturbation of $f$, using Theorem \ref{ArBC},  we can also suppose that $H(p,f)=M$. Using Lemma \ref{ex}, we can consider the following periodic linear system with symplectic transitions: $(\Sigma, \ f, \ TM,\ Df)$, where $\Sigma$ is the set formed by all the hyperbolic periodic points of $f$ homoclinically related to $p$. By Lemma \ref{exp}, there exists a periodic point $p_1\in \Sigma$ and a diagonalizable symplectic perturbation $A$ of $Df$ along the orbit of $p_1$, such that 
\begin{equation}
\left|S(p,f)-S(p_1,A)\right|<\eps/2.
\label{eqqq}\end{equation}
Since this perturbation is arbitrary close to $Df$, we can use Franks Lemma (Lemma \ref{franks}) to perform a non-linear perturbation of $f$ to $f_1\in \mathcal{U}$ such that $p_1$ is a periodic point of $f_1$, and moreover $Df_1$ along the orbit of $p_1$ coincides with $A$. In particular, $p_1$ is a diagonalizable periodic point of $f_1$, and $S(p_1,f_1)=S(p_1,A)$. As before, after a perturbation we can suppose that $H(p_1,f_1)$ coincides with $M$. And thus, we can consider the following periodic linear system: $(\Sigma_1, \ f_1, \ TM,\ Df_1)$, where $\Sigma_1$  is now the set formed by all the hyperbolic periodic points of $f_1$ homoclinically related to $p_1$.

Since $f_1\in \mathcal{U}$,  it has no dominated splitting. In particular, the above periodic linear system  has either no dominated splitting. Which implies, by Proposition \ref{every.rank}, that for any $1\leq i\leq 2n-1$ there is an arbitrary linear perturbation of $Df_1$ having complex eigenvalue of rank $(i,\ i+1)$. Therefore, this linear periodic system is in the hypothesis of Proposition \ref{blableblu}, which implies the existence of a periodic point $\tilde{p}\in \Sigma_1$, and a symplectic perturbation $\tilde{A}$ of $Df_1$ along the orbit of $\tilde{p}$ such that:

\begin{equation} 
\left|n\chi^+_{min}(\tilde{p},\tilde{A})-S(p_1,f_1)\right|<\eps/2.
\label{eqqqqq}\end{equation}

Once again, we can use Franks Lemma (Lemma \ref{franks}) to perform a local perturbation of $f_1$ to $g\in 
\mathcal{U}$ along the orbit of $\tilde{p}$ such that: $\tilde{p}$ still is a periodic point of $g$, and $Dg$ coincides with $\tilde{A}$ in the orbit of $\tilde{p}$. In particular $\tilde{p}$ is a hyperbolic periodic point of $g$, and moreover, from equations (\ref{eqqq}),  and $(\ref{eqqqqq})$ we have  $\left|n\chi^+_{min}(\tilde{p},g)-S(p,f)\right|<\eps$.

$\hfill\square$

\section{Abundance of full homoclinic tangencies: proof of Proposition \ref{fulltangency}.}

In this section we will prove Proposition \ref{fulltangency}. That is, we prove that in the lack of partial hyperbolicity any hyperbolic periodic point exhibits a  full homoclinic tangency, after a perturbation. Recall a hyperbolic periodic point $p$ of a $C^1-$diffeomorphism $f$ exhibits a homoclinic tangency if there is a non-transversal intersection $q$ between its stable and unstable manifold, $W^s(p,f)$ and $W^{u}(p,f)$. We define
$c_T(q) = dim (T_qW^u(p,f)\cap T_qW^s(p,f))$
as  the {\it codimension of the tangency at $q$}. According with such definition we can say we have a {\it full homoclinic tangency} when we have a tangency having the maximal possible codimension. In particular if $f\in \dw$, then a hyperbolic periodic point $p$ has a full homoclinic tangency if exhibits a homoclinic tangency with codimension $n$.  

We would like to point out that recently Barrientos and Raibekas proved in \cite{BR} the existence of robust homoclinic tangencies with large codimension. Their result also works for symplectic diffeomorphisms. 

However, to prove Proposition \ref{fulltangency}  we will use some ideas developed in \cite{CH}. But, it is worth to point out, that in our setting, we do not get robustness. 

Before we start the proof we need some definitions and results. 

Considering $(M^{2n},\omega)$ a symplectic manifold, we say that a sub manifold $N\subset M$ is Lagrangian if $dim\ N=n$ and moreover if $\omega(x)|T_xN\times T_xN=0$ for any $x\in N$. In particular, if $p$ is a hyperbolic periodic point of $f\in \dw$, then its stable and unstable manifolds are Lagrangian sub manifolds of $M$. 

In the proof bellow we will need to connect stable and unstable manifolds of periodic points. To do this we will use the famous connecting lemma of Hayashi \cite{H}. The next result is a symplectic version of it, which was proved by Wen and Xia \cite{XW}. 

\begin{theorem}[$C^1$-connecting lemma]
Let $f \in \dw$ and $p_1,\, p_2$ hyperbolic periodic points of $f$, such that there exist sequences  $y_n\in M$ and positive integers 
$k_n$ such that:
\begin{itemize}
\item $y_n\rightarrow y \in W_{loc}^u(p_1, f))$, $y\neq p_1$; and
\item $f^{k_n}(y_n)\rightarrow x \in W_{loc}^s(p_2, f))$, $x\neq p_2$.
\end{itemize}
Then, there exists a $C^1$ symplectic diffeomorphism  $g$,   $C^1-$close to $f$, such that $W^u(p_1,g)$ and $W^s(p_2,g)$ have a non empty intersection close to  $y$. \label{connecting lemma}\end{theorem}

\vspace{0,3cm}

{\it Proof of Proposition \ref{fulltangency}: }
Let $\mathcal{U}\subset \dw$ be a small neighborhood of $f$ and $p$ be a hyperbolic periodic point of $f$, as given in the proposition. Now, taking $\mathcal{U}$ smaller, if necessary, we can assume that for any $g\in \mathcal{U}$: 
\begin{itemize}
\item the analytic continuation of $p$ is defined for $g,$ and

\item $g$ is not partially hyperbolic.
\end{itemize}

Now, we can use Proposition 5.3 in \cite{HT} to find a diffeomorphism $f_1\in \mathcal{U}$ having a periodic point $\tilde{p}$ such that $Df_1^{\tau(\tilde{p},f_1)}=Id$.

Let $\gamma>0$ be an arbitrary small positive real number. By the same arguments as in the proof of Theorem A in \cite{CH} we can perform a local perturbation of $f_1$ to $f_2\in \mathcal{U}$ such that $\tilde{p}$ is now a hyperbolic periodic point of $f_2$ and moreover $W^s_{loc}(\tilde{p},f_2)$ is a Lagrangian sub manifold of $M$ $\gamma-C^1$ close to the Lagrangian sub manifold  $W^u_{loc}(\tilde{p},f_2)$. Let $\tilde{\mathcal{U}}\subset \mathcal{U}$ be a small neighborhood of $f_2$ such that the analytic  continuation of $\tilde{p}$ is defined in $\tilde{\mathcal{U}}$ and moreover the local stable and unstable manifolds of $\tilde{p}(g)$ still are sub manifolds $\gamma-C^1$ close, for any $g\in \tilde{\mathcal{U}}$. 

The next step of the proof is similar to the proof of Lemma 5.2 in \cite{CH}. In words, we want to perturb $f_2$ in order to find a transversal homoclinic point $q$ of $p$ near $\tilde{p}$, such that the stable and unstable manifolds of $p$ at $q$ are locally  $C^1-$close to the local stable and unstable manifolds of $\tilde{p}$, respectively.  

 First,  we  use Theorem \ref{ArBC} to perturb $f_2$ to a transitive symplectic diffeomorphism $f_3\in \tilde{\mathcal{U}}$. Hence, we can use Theorem \ref{connecting lemma} to perturb $f_3$ to $f_4$ such that  there is an intersection between $W^s(p(f_4))$ and $W^u(\tilde{p}(f_4))$. After a perturbation, if necessary, we can suppose this intersection is transversal. By transversality and  continuity of stable and unstable manifolds of hyperbolic periodic points, this intersection is robust. Hence, we can use again Theorem \ref{ArBC} and Theorem \ref{connecting lemma} to perturb $f_4$ to  a diffeomorphism $\tilde{g}\in \tilde{\mathcal{U}}$ such that for $\tilde{g}$  there is also an intersection between $W^u(p(\tilde{g}))$ and $W^s(\tilde{p}(\tilde{g}))$.  Now, since $\tilde{p}(\tilde{g})$ and $p(\tilde{g})$ have the same indices (i.e. their stable manifolds have same dimension), the existence of such transversal intersections between the stable and unstable manifolds of $\tilde{p}(\tilde{g})$ and $p(\tilde{g})$, together with the famous Lambda-Lemma, allow us to find a point $q\in W^s(p(\tilde{g}))\cap W^u(p(\tilde{g}))$, such that: 
 \begin{itemize}
\item  $T_qW^{s}(p(\tilde{g}))$ is $\gamma-$close to $E^{s}(\tilde{p}(\tilde{g}))=T_{\tilde{p}}W^{s}(\tilde{p}(\tilde{g}))$, and  

\item  $T_qW^{u}(p(\tilde{g}))$ is $\gamma-$close to $E^{u}(\tilde{p}(\tilde{g}))=T_{\tilde{p}}W^{u}(\tilde{p}(\tilde{g}))$.   
  \end{itemize}
And thus  $T_qW^{s}(p(\tilde{g}))$ is $2\gamma-$close to $T_qW^{u}(p(\tilde{g}))$.   

Now, we have the following lemma: 

\begin{lemma}
\label{isotropic2}
Suppose a symplectic vector space $V=E\oplus F$, where $E$ is a Lagrangian subspace.
For any $\delta>0$, there exists $\gamma>0$ such that if $W\subset V$ is a
Lagrangian  subspace $\gamma$-close to $E$, then there exists a symplectic
linear map $B$ on $V$ $\delta$-close to $Id$ such that $B(W)=E$ and
$B|F=Id_F$.
\end{lemma}
 
There is a proof of such lemma in Section 2 of \cite{CH}. 

Therefore, given $\delta>0$ we can choose $\gamma>0$ as in Lemma \ref{isotropic2}, and thus, since the stable and unstable manifolds of $p$ are Lagrangian Manifolds, there exists  a  symplectic linear map $B$ in $T_qM$, $\delta-C^1$ close to identity, such that $B(T_qW^{u}(p(\tilde{g})))=T_qW^{s}(p(\tilde{g}))$. 

Finally, we can use Franks Lemma, Lemma \ref{franks},  to perform a local perturbation in a small neighborhood of $\tilde{g}^{-1}(q)$, and find $g\in\tilde{\mathcal{U}}$ such that $q$ still is a homoclinic point of $p$ and moreover $Dg(\tilde{g}^{-1}(q))=B\circ D\tilde{g}(\tilde{g}^{-1}(q))$. Which implies that $T_qW^s(p(g))=T_qW^u(p(g))$, as we wanted. 

$\hfill\square$

\section{Snake perturbation product and proof of Proposition \ref{tech.prop}.}

In this section we prove the key technical result (Proposition \ref{tech.prop}) in the proof of Theorem \ref{hgeqs}. 
To do that we improve the snake perturbation created by Newhouse \cite{N}, and used by many authors since then. For instance, we found a way to make snake perturbations simultaneously in many directions.

In this section we will need to perform local $C^1-$perturbations  of symplectic diffeomorphisms. More precisely, we need to paste distinct maps, locally.  For general diffeomorphisms this is quite easy to do by using bump functions. However, since we are in the symplectic setting, we need to keep the symplectic structure after we paste the maps. Fortunately, we can do this, locally, by using generating functions. More precisely,   

\begin{lemma}[Lemma 3.9 in \cite{AM}]  If $f$ is a $C^k$-symplectic diffeomorphism on $M^{2n}$ ($k\geq 1$), then for any $x\in M$ and $g$ a local sympletic diffeomorphism ($C^k$-close to $f$ ) defined in a small neighborhood $U$ of $x$, there exists a $C^k$-symplectic diffeomorphism $h$ ($C^k$-close to $f$) and some neighborhood $V\subset  U$ of $x$ satisfying $h|V = g$ and $g|U^c=f$. 
\label{past2}\end{lemma}

For more details about generating functions see also \cite{AM}.

\vspace{0,2cm}
{\it Proof of Proposition \ref{tech.prop}: } Let $f\in \dw$ and $p$ be a hyperbolic periodic point of $f$ exhibiting a full homoclinic tangency at $q$, as in the hypothesis of the proposition.  

We can take a $C^1$ neighborhood $\mathcal{U}\subset \dw$ of $f$ small enough, such that there exists a common Lipschitz constant, say $\lambda>0$,  for every diffeomorphism in $\mathcal{U}$. To simplify notation we will use in this proof $\tau=\tau(p,f)$. 

Now, we will linearize the diffeomorphism in a small neighborhood of $p$.  
After a local perturbation, by transversality,  we can assume that $q$ is a transversal homoclinic point but having the spaces $T_qW^s(p)$ and $T_qW^u(p)$ arbitrary close to each other, since the stable and unstable manifolds of a hyperbolic periodic point varies continuously in compact parts. We use now Pasting Lemma, Lemma \ref{pasting lemma}, to perform a $C^1-$perturbation of $f$ to $f_1\in \mathcal{U}$ such that $p$ still is a hyperbolic periodic point of $f_1$ and there is a small neighborhood $V$ of $p$ such that $f^{\tau}_1|V=Df^{\tau}(p)$ (in local coordinates) and $q\not\in V$. By continuity of the stable and unstable manifolds of $p$, $f_1$ still exhibits a transversal homoclinic point $\tilde{q}$ of $p$ close to $q$, where the tangent spaces to $W^s(p)$ and $W^u(p)$ at $\tilde{q}$ still are near.  
Hence, after another local $C^1$ perturbation of $f_1$ in a small neighborhood of $\tilde{q}$, as done in the proof of Proposition \ref{fulltangency},  we can recover the homoclinic tangency to $f_1$ at $\tilde{q}$ as in the hypothesis of the propositon.

 Replacing $\tilde{q}$ by an iterate we can suppose that $\tilde{q}\in V$ and $f_1^{-1}(\tilde{q})\not\in V$. 
Thus,  we can choose a small open set $U\subset V$  containing  $\tilde{q}$ such that $f_1^{-1}(U)\cap V=\emptyset$. Since $f^{\tau}_1$ is linear on $V$ (in local coordinates)  we can choose a local coordinate $\psi$ of $U$ such that  $\psi(W^s_{loc}(p)\cap U)\subset R^n \times \{0\}^n$, and such that $\psi(\tilde{q})=0$. Let $\Xi^u$ be the connected component of $W^u(p)\cap U$ containing $\tilde{q}$. After shrinking $U$, if necessary, we have that $\psi(\Xi^u)$ is a graphic of a map  $\xi: \R^n\rightarrow \R^n$, with zero derivative at origin since $\tilde{q}$ is a homoclinic tangency. Now, since the stable and unstable manifold of $p$ are Lagrangian sub manifolds  of $M$, the map $\pi(x,y)=(x,y-\xi(x))$ ($x,y\in \R^n$) defined in $\psi(U)$ preserves the symplectic structure and moreover is arbitrary $C^1-$close to identity in a small neighborhood of the origin.  Thus, we can use Lemma \ref{past2} to perform a local $C^1$ perturbation of $f_1$ to find a diffeomorphism $f_2\in \mathcal{U}$ such that $f_2$ is equal to $\psi^{-1}\circ \pi\circ \psi\circ f_1$ in a small neighborhood of $f_1^{-1}(\tilde{q})$. In particular, $f_2$ satisfies: 
\begin{itemize}
\item[1-] $f_2(x)=f_1(x)$ for any  $x$ in the complement of $f_1^{-1}(U)$;

\item[2-] There are $a>0$ and a disk $\tilde{D}=[-a,a]^n\times \{0\}^n\subset \R^{2n}$, such that $D^s=\psi^{-1}(\tilde{D})\subset W^s_{loc}(p)\cap W^u(p)$.  
\end{itemize}    

In words, item 2 says that $f_2$ has a disc $D^s$ of homoclinic points of $p$. 

Since $D^s$ is a small disc of homoclinic points of $p$ around $\tilde{q}$, let $T$ be a positive integer such that $D^u=f_2^{-T}(D^s)\subset V$. Note that $D^u$ is a disc inside $W^u_{loc}(p)$.

The next step of the proof is to perform a local $C^1$ perturbation of $f_2$ in order to create arbitrary finitely many transversal  homoclinic points of $p$ inside $U$. This kind of perturbation is known by snake perturbation, and was first introduced by Newhouse \cite{N}. However, here we repeat his idea in many directions, which we call by {\it snake perturbation product.}

For that,  we fix $\delta>0$ arbitrary small, such that any diffeomorphism $C^1$ $\delta-$close  to $f$ must be in $\mathcal{U}$. Let $N$ be a large odd natural number and $\tilde{A}>0$ be a small real number. Then we consider the following symplectic map defined in $\R^{2n}$:  
$$
\Theta_{\tilde{A},N}(x_1,\ldots,x_n,y_1,\ldots, y_n)=(x_1,\ldots,x_n, y_1+\tilde{A}\sin \displaystyle\frac{x_1\pi N}{2a}, \ldots, y_n+\tilde{A}\sin \displaystyle\frac{x_n\pi N}{2a}). 
$$
It is not difficult to see that $\Theta(\tilde{D})\cap \tilde{D}$ is a finitely set containing $N^n$ points. Moreover note that $\Theta_{\tilde{A},N}$ is $C^1$ $\displaystyle\frac{\tilde{A}\pi N}{2a}-$close to the identity map. Hence, considering $A=\displaystyle\frac{2a\alpha\delta}{\pi N}$, where $\alpha$ is a constant depending on the local coordinate $\psi$, we have that the map $\psi^{-1}\circ \Theta_{A,N} \circ\psi$ defined in $U$ is $C^1$ $\delta-$close to identity map restricted to $U$. Thus, we can use Lemma \ref{past2} to make a local perturbation of $f_2$ to $f_{3,N}\in \mathcal{U}$ such that: 
\begin{itemize}
\item $f_{3,N}(x)=f_2(x)$ for any  $x$ in the complement of $f_2^{-1}(U)$;

\item $f_{3,N}$ is equal to $\psi^{-1}\circ \Theta_{A,N}\circ \psi\circ f_2$  in an open set $\tilde{U}\subset U$ which contains the disc $D^s$.   
\end{itemize}    

Hence, by construction of  $\Theta_{A,N}$, we remark that this $C^1$ local perturbation destroys the disc  of homoclinic points $D^s$ of $p$ by creating $N^n$ transversal homoclinic points. We will denote $J=\psi^{-1}\circ \Theta_{A,N}\circ \psi(D^s)$, which is a disc in $W^u(p)$ containing these $N^n$ homoclinic points. It is also important to remark that by properties of the unstable and stable manifold it is true that  $f_{3,N}^{-T}(J)=D^u\subset W_{loc}^u(p)$, and $D^s$ sitll belongs to $W^s_{loc}(p)$, since $f_{3,N}$ was obtained by a $C^1$ local perturbation. 

We will use these transversal homoclinic points to find a hyperbolic set $\Lambda_p$ satisfying the thesis of proposition. For that, we consider  a local coordinate map $\phi$ on $V$ such that  $\phi(p)=0$ and $\phi\circ f^{\tau}_{3,N}\circ\phi^{-1}$ is linear. Also, we consider a decomposition in $\R^{2n}$ induced by the hyperbolic splitting  $E^s(p)\oplus E^u(p)$ existent in $T_{p}M$ for $Df^{\tau}_{3,N}(p)$, using the local coordinate $\phi$. To simplify notation we continue denoting this splitting by  $E^s(p)\oplus E^u(p)$.  Using such decomposition  we define $\Gamma_l=B^s_l\times D^u$ and $b_l^s(x)=B_l^s\times \{x\}$ for any $x\in \phi(D^u)$. Here $B^s_l$ is the closed ball with centre zero and radius $l$ inside $E^s(p)$.  
We will fix  $l_D$ such that $B^s_{l_D}$ contains the disc $D^s$.

Since $f_{3,N}\in\mathcal{U}$, $\lambda$ is a Lipschitz constant to $f_{3,N}$. Hence, if $l_A=\displaystyle\frac{A\beta}{4\lambda^T}$  then  
 $f_{3,N}^T(\phi^{-1}(\Gamma_{l_A}))$ is a neighborhood of $J$ contained in a $A/4-C^1$ neighborhood of $J$. Here, $\beta$ is a constant which depends only on $\phi$.   

We choose now the smallest positive integer $t$  such that:
\begin{align}
&\hspace{-1cm}  f_{3,N}^{-\tau t}(\phi^{-1}(\Gamma_{l_A}))\cap V \ \text{ is }\  A/4-C^1 \ \text{ close to }\  W^s_{loc}(p), \text{ and }
\label{t.1}\\
 &\hspace{-1cm} \text{ the projection of }\ \phi\circ f_{3,N}^{-\tau t}\circ \phi^{-1}(b^s_{l_A}(x)) \ \text{ in }\  E^s(p) \ \text{ contains the subset } \ B^s_{l_D}\cap \phi(U).\label{t.2}  
\end{align}

The next lemma says there is a relation between $t$ and $A$.

\begin{lemma}
 \label{afirma}
For $A$ and $t$ defined as before, there exists  a positive integer $K$, which is 
independent of $A$ (in particular of $N$), such that $$K^{-1} \min \{ \|Df^{t\tau}(p)|E^u\|^{-1},\,
\|Df^{-t\tau}(p)|E^s\|^{-1}\}\leq A\leq K \max \{ \|Df^{-t\tau}(p)|E^u\|,\,
\|Df^{t\tau}(p)|E^s\|\}.$$\end{lemma}

Proof: By construction of $f_{3,N}$, note that  $\phi\circ f^{\tau}_{3,N}\circ \phi^{-1}$ is linear and equal to $D\phi(p)\ Df^{\tau}(p)\ D\phi^{-1}(0)$ in $\phi(V)$. 
In particular, the action of $f^{\tau}_{3,N}$ on $V$ does not depend on $N$. Also, recall that  $D^u$ and $l_{D^s}$ do not depend on $N$, too.

Hence, by the above remarks and the choice of $l_A$ there exists constants $K_1$ and $K_2$, which do not depend on $N$ (in particular do not depend on $A$) such that for any $0\leq i\leq t$ and $x\in f_{3,N}^{-i\tau}(f^{i\tau}_{3,N}(\phi^{-1}(\Gamma_{l_A}))\cap V)$:

\begin{equation}
K^{-1}_1\|Df^{i\tau}(p)|E^u\|^{-1}\leq d(f^{-i\tau}_{3,N}(x),W^s_{loc}(p))\leq K_1 \|Df^{-i\tau}(p)|E^u\|, \text{ and } \label{desigualdade 1}\end{equation} 
\begin{equation}
K^{-1}_2 \ A\ \|Df^{i\tau}(p)|E^s\|^{-1}\leq d(f^{-i\tau}_{3,N}(x),W^u_{loc}(p))\leq K_2\  A\ \|Df^{-i\tau}(p)|E^s\|. \label{desigualdade 2}\end{equation}

Thus, using the above inequalities for a point $x\in f_{3,N}^{t\tau}(f_{3,N}^{-\tau t}(\phi(\Gamma_{l_A}))\cap V)$ together with (\ref{t.1}) and (\ref{t.2}) we obtain: 
$$A\geq \min\{4K^{-1}_1\|Df^{t\tau}(p)|E^u\|^{-1}, \; l_{D^s} K^{-1}_2 \|Df^{-t\tau}(p)|E^s\|^{-1} \}
$$

On the other hand, by  choice of $t$, there exists $x\in \phi(\Gamma_{l_A})$ such that $ f^{i\tau}_{3,N}(x)\in V$ and
\begin{itemize}
\item either $d(f^{-(t+1)\tau}_{3,N}(x),W^s_{loc}(p))\geq\displaystyle\frac{A}{4}$, 

\item or $d(f^{-(t+1)\tau}_{3,N}(x),W^u_{loc}(p))\leq l_{D^s}.$
\end{itemize}
Now, these two above inequities together  with the inequalities (\ref{desigualdade 1}) and (\ref{desigualdade 2}) imply that 
$$
A\leq \max\{4 K_1 \|Df^{-(t-1)\tau}(p)|E^u\|, \; l_{D^s} K_2 \|Df^{(t-1)\tau}(p)|E^s\| \}
$$

Therefore, the proof is finished taking $$K=\max\{4 K_1 \|Df^{\tau}(p)|E^u\|, \; l_{D^s} K_2 \|Df^{-\tau}(p)|E^s\|,\; K_1/4, \; K_2/l_{D^s}\}$$

$\hfill\square$ 
 
Finally, using $t$ we will construct a hyperbolic set close to the disc $D^s$ and using the previous lemma we will estimate the topological entropy of such hyperbolic set by means of the Lyapunov exponents of $p$. 

For that we define  $\tilde{R}=\phi(D^s)\times \phi(f_{3,N}^{-\tau t}(D^u))$ and $R=\phi^{-1}(\tilde{R})$.  Thus, the choice of $t$ implies that $f_{3,N}^{T+\tau t}(R)$ intersects transversally $R$ in $N^n$ disjoint connected components. Then, it is known that the maximal invariant set in $R$ is a basic hyperbolic set $\tilde{\Lambda}_{p,N}$ which is conjugated to the shift acting in the space of sequences of $N^n$ symbols. In particular $h_{top}( f_{3,N}^{-T-\tau t}|\tilde{\Lambda}_{p,N})=\log N^n$. Moreover, if we define $\Lambda_{p,N}=\cup_{i=0}^{T+t\tau-1}  f_{3,N}^{i}(\tilde{\Lambda}_{p,N})$, we have that 
\begin{equation}
h_{top}( f_{3,N}|\Lambda_{p,N})=\displaystyle\frac{1}{T+\tau t} \log N^n. 
\label{porrada}\end{equation}

It is important to remark now, that  a directly consequence of Lemma \ref{afirma} is that $t$ goes to infinity when $N$ goes, since  $A=\displaystyle\frac{2a\alpha\delta}{\pi N}$. Hence, taking  $\eps>0$ arbitrary small we can use this information together with Lemma \ref{afirma} to choose a large  positive integer $N$, such that
$$
\displaystyle\frac{1}{T+\tau t}\log N> \min\left\{\frac{1}{T+\tau t}\log\|Df^{-\tau t}(p)|E^u\|^{-1},\, \frac{1}{T+\tau t}\log\|Df^{\tau t}(p)|E^s\|^{-1} \right\} -\frac{\eps}{2}.
$$
Then, when $t$ goes to infinity (which happens when $N$ goes) the minimum of the right side of the above inequality converges to $\chi^+_{min}(p)$ (the smallest positive Lyapunov exponent of $p$),
by definition. Thus, there exists a large positive integer $N_0$ such that
\begin{equation}
\displaystyle\frac{1}{T+\tau t}\log N_0> n\ \chi^+_{min}(p)-\eps.
\label{porradaa}\end{equation}
Therefore, taking $g=f_{3,N_0}$ and $\Lambda_p=\Lambda_{p,N_0}$ we can use (\ref{porrada}) and (\ref{porradaa}) to conclude:

$$
h_{top}(g|\Lambda_p)> n\ \chi^+_{min}(p,g)-\eps,$$
which finishes the proof. 

$\hfill\square$

\bigskip

\flushleft

{\bf Thiago Catalan} (tcatalan\@@famat.ufu.br)\\
Faculdade de Matem\'{a}tica, FAMAT/UFU \\
Av. Jo\~ao Naves de Avila, 2121\\
38.408-100, Uberl\^andia,MG, Brazil

\end{document}